\documentclass[11pt]{article}
\usepackage{amsfonts,amsmath,amssymb,amsthm, amscd}
\usepackage{graphicx}

\numberwithin{equation}{section}

\newtheorem{theorem}{Theorem}[section]

\newtheorem{lemma}[theorem]{Lemma}
\newtheorem{cor}[theorem]{Corollary}

\theoremstyle{definition}

\newtheorem{remark}[theorem]{Remark}

\def\<{{\langle}}
\def\>{{\rangle}}

\def\Z{\mathbb Z}

\DeclareMathOperator\vg{vg}

\begin{document}

\title{Virtual Genus of Satellite Links}

\author{Daniel S. Silver \and Susan G. Williams\thanks{The
work of both authors was partially supported by grants \#245671 and \#245615 from the Simons
Foundation.} \\ {\em
{\small Department of Mathematics and Statistics, University of South Alabama}}}

\maketitle 

\begin{abstract} The virtual genus of a virtual satellite link is equal to that of its companion. \\


\end{abstract}

\section{Introduction} 

The notion of virtual knots and links was introduced by L. Kauffman \cite{kauf99}. It is a nontrivial extension of the classical theory. Virtual links can be defined as link diagrams in the plane with ``virtual crossings" as well as crossings of the usual kind and an extended set of Reidemeister moves, or as combinatorial Gauss diagrams.
It is shown in \cite{kamkam00} that virtual links correspond bijectively to abstract link diagrams, introduced by N. Kamada in \cite{kam97}.

Alternatively, a virtual link $\ell$ can be defined as an equivalence class of link diagrams ${\cal D}$ in a surface $S$. The surface is required to be closed and orientable; it need not be connected, but we require that each component contain at least one link component. The equivalence relation is generated by Reidemister  moves on ${\cal D}$,  orientation-preserving homeomorphisms of $S$  and adding or deleting hollow 1-handles in the complement of the diagram. 

Adding a handle ({\it stabilization}) is the surgery that removes two open disks disjoint from ${\cal D}$, and then joins the resulting boundary components by an annulus. Deleting a handle ({\it destabilization}) is the surgery that removes the interior of a neighborhood of a simple closed curve  that misses ${\cal D}$, and then attaches a pair of disks to the resulting boundary. Destabilization might produce a  diagram for the link in a surface that has smaller (total) genus than $S$. 

Following \cite{dk}, we define the {\it virtual genus} of $\ell$, denoted here by $\vg(\ell)$, to be the minimal genus of a surface that contains a diagram representing the link. The virtual genus was first studied in \cite{kam97}, where it is called {\it supporting genus}.  Methods for estimating virtual genus are found in \cite{dk}, \cite{kauf}, \cite{csw}. 

By \cite{kauf99}, \cite{cks}, a virtual link can also be regarded as an equivalence class of embedded links in thickened surfaces. The equivalence relation is generated by isotopy as well as stabilization/destabilization. Destabilization in this context consists of surgery along an embedded annulus $A$ that is {\it vertical} in the sense that $A = p^{-1}(p(A))$, where $p$ is first-coordinate projection $S \times I \to S$. The reverse operation, stabilization, is a parametrized connected-sum operation with a thickened torus.\\

Lemma 3.4 of \cite{wald68} implies the following.
\begin{lemma}\label{vertical} A properly embedded annulus $(A; \partial_1A, \partial_0 A) \subset (S\times I; S\times \{1\}, S\times \{0\})$ is isotopic to a vertical annulus provided $\partial_iA$ is essential in $S \times \{i\},\ i=0, 1$.
\end{lemma}

Theorem 1 of \cite{kup03} implies that if $\vg(\ell)$ is less than the genus of $S$, then after isotopy a vertical annulus $A \subset S \times I \setminus \ell$ can be found such that surgery along it produces an embedding of the link in a surface of strictly smaller genus than that of $S$. Note that surgery on such an annulus will reduce genus if and only if each of its boundary components represents a nontrivial element of $H_1(S; \Z)$. \\



Assume that $\ell = \ell_1 \cup \cdots \cup \ell_d \subset S \times I$ is a link in a thickened surface. Let $N = N_1 \cup \cdots \cup N_d$ be a regular neighborhood of $\ell$ with boundary $\partial N  = \partial N_1 \cup \cdots \partial N_d$ consisting of mutually disjoint tori. For each $i = 1, \ldots, d$, let $\tilde \ell_i \subset {\rm int}\ N_i$ be a link that is not contained in any 3-ball neighborhood in $N_i$. Then $\tilde \ell = \tilde \ell_1 \cup \cdots \cup \tilde \ell_d$ is a {\it satellite link} of $\ell$ with {\it companion} $\ell$. \\

It is clear that the virtual genus of $\tilde \ell$ is not greater than that of $\ell$. 
The following theorem asserts the virtual genus of the links are in fact equal. 

\begin{theorem} \label{main} If $\tilde \ell$ is any satellite link with companion $\ell$, then the virtual genus of $\tilde \ell$ is equal to that of $\ell$.
 \end{theorem}
 
 A virtual link $\ell$ is {\it classical} if $\vg(\ell)=0$ (equivalently, if it can be represented by a planar diagram).
 
 \begin{cor} If a satellite virtual link is classical, then its companion is classical. 
\end{cor}   
 
\begin{remark} 

The main result of \cite{kup03} is that every virtual knot has a unique representative $\ell \subset S \times I$ for which the genus of $S$ is equal to $\vg(\ell)$ and the number of components of $S$ is maximal.  (Uniqueness is up to isotopy and orientation-preserving self-homeomorphism of $(S\times I, S\times\{1\}, S\times\{0\})$.)  Let  $\tilde \ell \subset S\times I$ be a satellite with companion $\ell$, and assume that $S$ has no genus-0 components and no embedded 2-sphere in $N_i$ separates $\tilde \ell_i$. We see easily from the proof of Theorem \ref{main} that $S$ has both minimal genus and maximal number of components for $\tilde \ell$.  Hence topological invariants of $S\times I \setminus \tilde \ell$ (e.g. fundamental group, homology groups of abelian covers) are also invariants of the virtual link $\tilde \ell$.

 \end{remark} 

The authors are grateful to J. Scott Carter and Seiichi Kamada for helpful suggestions.
 
 \section{Proof of Theorem \ref{main}.}
 
Consider a link $\ell = \ell_1\cup \cdots \cup \ell_d  \subset S \times I$ such that the genus of $S$ is equal to 
 $\vg(\ell)$. Let $\tilde \ell = \tilde \ell_1 \cup \cdots \cup \ell_d\subset S \times \tilde I$ be any satellite link with companion $\ell$, as above. 
 
Suppose that some embedded 2-sphere in $S \times I$ separates $\ell$ into  nonempty sublinks. Since one of the sublinks must be contained in a 3-ball, it suffices to prove Theorem \ref{main} for the other sublink.
By an induction argument, we can assume without loss of generality that no 2-sphere in $S \times I$ separates $\ell$.  

Similarly, we may assume that no embedded 2-sphere or properly embedded annulus in $N_i$ separates $\tilde \ell_i$.  Otherwise, at least one of the two sublinks of $\tilde \ell_i$ is not contained in any 3-ball, and it suffices  prove the result for the link obtained by deleting the other sublink.

Assume that $\vg(\tilde \ell) < \vg(\ell)$. By the argument of \cite{kup03}, there exists a  vertical annulus
$$(A; \partial_0 A, \partial_1 A) \subset (S\times I\setminus \tilde \ell; S \times \{0\}, S \times \{1\})$$ such that $\partial_i A$ represents a nontrivial element of $H_1(S\times \{i\}; \Z)$ for $i = 0, 1$. We will derive a contradiction. 

Deform $A$ so that it meets $\partial N$  transversely. Regard the intersection as a closed 1-submanifold of $A$. 

If some component of $A \cap \partial N$ is null-homotopic in $A$, then let
$C$ be such a component that is innermost in the sense that no other component of $A \cap \partial N$ is contained in the 2-disk $D \subset A$ bounded by $C$. Let $\partial N_i$ be the component of $\partial N$ that contains $C$. 

The assumption that $\tilde \ell_i$ is not contained in a 3-ball in $N_i$, implies that $\partial N_i$ is incompressible in $N_i \setminus \tilde \ell_i$. Hence $C$ also bounds a 2-disk $D'$ in $\partial N_i$. The union $D \cup_C D'$ is an embedded 2-sphere bounding a ball in $N_i$ that does not meet $\tilde \ell_i$, by our assumption.  We can deform $A$ in this ball to remove the circle $C$ of intersection with $\partial N_i$. By repeating this procedure, we can assume that every component of $A \cap \partial N$ is essential in $A$. 

Now $A \cap \partial N$ consists of finitely many pairwise disjoint simple closed curves that divide $A$ into successive annular regions $A_1, B_1, \ldots, A_n, B_n, A_{n+1}$ such that each $A_j$ is contained in the exterior $X= S \times I \setminus {\rm int}\ N$ while each $B_j$ is contained in some $N_{i_j}\setminus \tilde \ell_{i_j}$. 

Consider any annular region $B_j$.  Its boundary components $\partial_\pm B_j$ are homologous in $N_{i_j}\setminus \tilde \ell_{i_j}$, and each is
homologous in $S\times I\setminus \ell$ to $\partial_1 A$. 
They separate $\partial N_{i_j}$ into two annuli
$A', A''$. Both $B_j\cup_\partial A'$ and $B_j\cup_\partial A''$ are tori; by our assumption, one of them contains the link $\tilde \ell_{i_j}$ while the other bounds a solid torus in $S \times I \setminus \tilde \ell$. The latter can be used to deform $A$ and remove the points $\partial_\pm B_j$ of intersection with $\partial N_{i_j}$. By repeating the procedure, we can deform $A$ away from $\partial N$.

The link $\ell$ is the core of $N$. Since $A \cap N$ is empty, we can use the vertical annulus $A$ to reduce the genus of $S$, contradicting the assumption that the genus of $S$ is equal to $\vg(\ell)$.


\end{document}